\newif
\ifconfver
\confvertrue        
\ifconfver
\documentclass[letterpaper, 10 pt, conference]{ieeeconf}  

\usepackage{xspace,empheq,fancybox,amssymb,amsfonts,graphicx,epstopdf,epsfig,syntonly,times,subfigure,caption} 
\usepackage{psfrag,color,bm,array,tabularx}
\usepackage{cite,url,footnote,xspace,syntonly,algorithmic}
\usepackage{verbatim,multirow}
\usepackage[linesnumbered,ruled,vlined]{algorithm2e}
\usepackage{textcomp}
\usepackage{xcolor}
\usepackage{mathtools}
\usepackage{bbm, booktabs}
\usepackage{amsmath, amsthm}           
\newcolumntype{M}[1]{>{\centering\arraybackslash}m{#1}}
\newcolumntype{N}{@{}m{0pt}@{}}
\DeclareMathOperator*{\argmax}{\arg\max}

\def\BibTeX{{\rm B\kern-.05em{\sc i\kern-.025em b}\kern-.08em
		T\kern-.1667em\lower.7ex\hbox{E}\kern-.125emX}}

\newtheorem{remark}{\bfseries Remark}
\input{mysymbol.sty}
\allowdisplaybreaks

\newcommand{\KB}{\color{black}{}}

\IEEEoverridecommandlockouts
\begin{document}
	

\title{Model-free Learning for Risk-constrained Linear Quadratic Regulator \\ with Structured Feedback in Networked Systems}
\author{\thanks{This work has been partially supported by NSF Grants 1802319,  1952193, and 2130706.} 
Kyung-bin Kwon\thanks{K. Kwon and H. Zhu are with the Department of Electrical \& Computer Engineering, The University of Texas at Austin, 2501 Speedway, Austin, TX, 78712, USA; Emails: {\{kwon8908kr, haozhu\}{@}utexas.edu}.}, 
Lintao Ye\thanks{L. Ye is with the School of Artificial Intelligence and Automation, Huazhong University of Science and Technology, 1037 Luoyu Road, Wuhan City, Hubei Province, 430074, China; Email: yelintao93@hust.edu.cn}, 
Vijay Gupta\thanks{V. Gupta is with the Department of Electrical Engineering, University of Notre Dame, 270 Fitzpatrick Hall, Nortre Dame, IN, 46556, USA; Email: vgupta2@nd.edu}, and Hao Zhu}
\maketitle

\begin{abstract}
We develop a model-free learning algorithm for the infinite-horizon linear quadratic regulator (LQR) problem. Specifically, (risk) constraints and structured feedback are considered, in order to reduce the state deviation while allowing for a sparse communication graph in practice. By reformulating the dual problem as a nonconvex-concave minimax problem, we adopt the gradient descent max-oracle (GDmax), and for model-free setting, the stochastic (S)GDmax using zero-order policy gradient. By bounding the Lipschitz and smoothness constants of the LQR cost using specifically defined sublevel sets, we can design the stepsize and related parameters to establish convergence to a stationary point (at a high probability). Numerical tests in a networked microgrid control problem have validated the convergence of our proposed SGDmax algorithm while demonstrating the effectiveness of risk constraints. The SGDmax algorithm has attained a satisfactory optimality gap compared to the classical LQR control, especially for the full feedback case. 
	
\end{abstract}


\section{Introduction}\label{sec:IN}

The linear quadratic regulator (LQR) problem is one of the most fundamental problems in optimal control theory \cite{LQR1,LQR2}. Recently, there is significant interest in model-free learning of the standard LQR problem using gradient-based approaches \cite{bu,malik}, with connection to the popular reinforcement learning (RL) methods. Nonetheless, model-free learning and convergence analysis for general LQR problems are still lacking such as constrained LQR and structured feedback design. 

Constraint functions have attracted recent interest for both LQ \cite{tsiamis,tsiamis2,zhao} and general RL problems \cite{paternain2019constrained,ding2020natural}. Constraints can increase the safety of the resultant policy while potentially improving the learning rates as a regularization. In particular, recent work \cite{tsiamis,zhao} has considered the mean-variance risk for the LQR problems, that can effectively mitigate the random state deviation from its mean due to noisy disturbance. More interestingly, \cite{tsiamis} has shown that this risk measure is equivalent to a quadratic constraint function that is similar to the LQR cost. In addition,  \cite{zhao} has developed a dual-ascent based double-loop algorithm by utilizing the global convergence of LQR learning \cite{bu,malik} for the inner-loop. Nonetheless, this double-loop procedure may be complicated to implement in practice because the inner-loop convergence is in a probabilistic way due to the stochastic gradient.  

Meanwhile, decentralized control problems \cite{shah2013cal,DEC,ye2021sample} arise in various real-world applications where sensors and actuators are distributed in a networked system. For example, it is very useful for power system control designs such as wide-area damping control \cite{dorfler2013sparse,chakrabortty2021wide} or networked microgird control \cite{DCM2,DLQR}. In decentralized LQR problems, a sparse communication graph leads to structured feedback gain, which has also been considered in recent gradient-based learning approaches \cite{bu,harvard}. In general, the stabilizable region of structured LQR is  disconnected with a complex geometry \cite{lavaeri}, and thus it is difficult to analyze. While gradient-based learning for structured LQR does not lead to global convergence as in the unstructured case \cite{bu,malik}, it is easy for implementation as the gradient can be simply performed over the non-zero entries \cite{bu,harvard}.

Our goal is to develop model-free learning algorithms for risk-constrained LQR problem under sparse feedback structure that arises in networked systems. The structured feedback is incorporated by considering the sparse non-zero entries only, and thus the gradient computation and updates can be performed without accounting for such structured constraint. Nonetheless, it leads to convergence to only a stationary point. As for the constraint function, it is similar to the LQR cost with the mean-variance risk as a special case as
shown by \cite{tsiamis,tsiamis2}. To deal with this constraint, we consider the dual problem which shares the stationary point (SP) with the minimax problem for the Lagrangian function. The resultant nonconvex-concave minimax reformulation motivates us to adopt Gradient-Descent max-oracle (GDmax) and the stochastic (S)GDmax algorithms in \cite{GDA} to solve the outer minimization problem via GD updates. More specifically, the SGDmax relies on the zero-order policy gradient (ZOPG) \cite{ZOO} which has bounded noise variance. 

Nonetheless, the key challenge in establishing the convergence results lies in the LQR cost function, which is shown to exhibit local-only Lipschitz and smoothness properties with location-dependent constants \cite{bu,malik}. To tackle this, we can introduce a compact sublevel set within which the upper bounds of Lipschitz and smoothness constants hold everywhere. Such analysis enables us to carefully design the stepsize and related parameters to establish the convergence to SP, while the convergence of SGDmax in a model-free setting can be attained with a high probability. Numerical results have validated the convergence of our algorithms and demonstrated the impact of having risk constraint and structured feedback in learning LQR policy. The SGDmax algorithm have attained satisfactory optimality gap compared to the classical LQR control, especially for the full feedback case.

The remainder of this paper is organized as follows. Sec.~\ref{sec:PF} formulates the infinite-horizon risk-constrained LQR with the structured feedback. Sec.~\ref{sec:DP} introduces the dual-related minimax reformulation and analyzes the convergence of the Gradient Descent with max-oracle (GDmax) algorithm. Sec.~\ref{sec:SG} extends it to model-free learning by developing the Stochastic (S)GDmax via zero-order policy gradient.  Sec.~\ref{sec:NT} presents the numerical results in a networked load frequency control (LFC) problem, while the paper is wrapped up in Sec.~\ref{sec:CN}.\\

\noindent \textbf{Notations:} Let $\|\cdot \|$ denotes the $L_2$-norm, $\nabla_{\ccalK}\ccalL$ the gradient of $\ccalL$ that admits  the structure defined in $\ccalK$, $\{X^j\}$ a sequence of $\{X^0, X^1, \ldots\}$, $\ccalP_{\ccalY}(\cdot)$ the projection onto the set $\ccalY$, and the operator $\otimes$ the Kronecker product of matrices. Last, $\mathbb E(\cdot)$ denotes the expectation while $\mathbb{P}(\cdot)$ the probability of an event. 


\section{Problem Formulation}\label{sec:PF}

We consider the infinite-horizon LQR problem for a linear time-invariant system, as given by
\begin{align}
	x_{t+1} = A x_t+B u_t+w_t, \quad t=0,1,\ldots \label{eq:dynamics}
\end{align}
with the state $x_t \in \mathbb{R}^n$, action $u_t \in \mathbb{R}^m$, and random noise $w_t \in \mathbb{R}^n$ that is uncorrelated across time. In addition, the model parameters $A \in \mathbb{R}^{n\times n}$ and  $B \in \mathbb{R}^{n\times m}$ can be unknown. The constrained LQR problem with structured feedback aims to  find an optimal linear feedback gain $K \in \mathbb{R}^{m \times n}$ for the control policy $u_t = -K x_t$ to:
\begin{align}
	\min_{K \in \mathcal{K}} \; R_0(K) &\! = \!\lim_{T \to \infty}\frac{1}{T}\mathbb{E} \sum_{t=0}^{T-1} [x^\top_t Q x_t + u^\top_t R u_t]  \label{eq:opt} \\
	\text{s.t.}\; R_i(K) &\!=\! \lim_{T \to \infty}\frac{1}{T} \mathbb{E}\sum_{t=0}^{T-1} [x^\top_t Q_i x_t + u^\top_t R_i u_t] \leq c_i,	\forall i \nonumber
\end{align}
{\KB where matrices $\{Q, R\}$ and $\{Q_i, R_i\}_{i\in \ccalI}$ are all positive semi-definite, with $\ccalI$ representing the set of the constraints.} 
The feasible set $\mathcal{K}$ enforces a structured policy, as
\begin{align}
	&\mathcal{K} = \{K: K_{a,b} = 0 \;\text{if and only if}\; (a,b) \notin \ccalE) \}
\end{align}
{\KB Here, the structure pattern $\ccalE$ is specified by the edges of a given communication or  information-exchange graph. Hence, the action for agent $a$, denoted as $u_{a,t}$, is determined as $u_{a,t} = -K_a x_{a,t}$, where $K_a$ is a row vector with only non-zero elements in $a$-th row of $K$ and $x_{a,t}$ is a sub-vector of $x_t$ according to $\ccalE$.} An example of the communication graph is illustrated in Fig.~\ref{fig:network}. {\KB The structured $\ccalK$ is  motivated by a multi-agent setting for networked control, where individual agents can  access partial feedback only depending on communication links. Notably, this structured constraint will lead to a complicated geometry of the feasible region \cite{bu,lavaeri}.} {\KB While the structured $\ccalK$ makes the analysis more difficult than the full feedback case, it does not increase the complexity of  computing the gradient as denoted by $\nabla_{\ccalK}$ later on. This is because one can represent the cost as a function of only non-zero entries in $K$ which can eliminate this structured constraint \cite{bu}. Accordingly,  the $\nabla_\ccalK$ operation needs no projection onto $\ccalK$, and can be thought of as the gradient for an unstructured $K$. Therefore, gradient-based methods are ideal for learning a structured policy.  }



As for the quadratic constraint in \eqref{eq:opt}, one can consider the mean-variance risk as a special instance,  represented by
\begin{align}
	R_c(K) \! =\!\! \lim_{T \to \infty}\!\frac{1}{T} \mathbb{E} \sum_{t=0}^{T-1} \left(x_t^\top Q x_t - \mathbb{E}[x_t^\top Q x_t \vert h_t]\right)^2 \leq \delta  \nonumber 
\end{align}
with the system trajectory $h_t := \{x_0, u_0, \ldots, x_{t-1}. u_{t-1}\}$ and a risk tolerance $\delta$. {\KB This risk measure limits the deviation from the expected cost given the past trajectory, and thus can mitigate extreme scenarios due to the uncertainty in the noisy dynamics.} Interestingly, under a finite  fourth-order moment of noise $w_t$, \cite{tsiamis,tsiamis2} has developed a tractable reformulation $R_c(K)$, as   
\begin{align}
	R_c(K) \!= \!\! \lim_{T \rightarrow \!\infty} \frac{1}{T} \mathbb{E} \sum_{t=0}^{T-1} \left(4x_t^\top QWQ x_t+4x_t^\top QM_3 \right) \!\leq\! \bar{\delta}\! \label{eq:objr2}
\end{align}
with $\bar{\delta} = \delta-m_4 + 4 \text{tr}\{(WQ)^2 \}$ and the (weighted) noise statistics given as
\begin{align}
	\bar{w}&=\mathbb{E}[w_t],\\
	W&=\mathbb{E}[(w_t-\bar{w})(w_t-\bar{w})^\top],\\ 
	M_3 &= \mathbb{E}[(w_t-\bar{w})(w_t-\bar{w})^\top Q (w_t-\bar{w})],\\
	m_4&=\mathbb{E}[(w_t-\bar{w})^\top Q (w_t-\bar{w}) - \text{tr}(WQ)]^2.
\end{align} 
With known noise statistics, this risk constraint shares the quadratic form in \eqref{eq:opt} with an additional linear term, which does not affect our proposed gradient-based learning. The ensuing section first develops  the deterministic algorithm for problem \eqref{eq:opt}, which can provide insights on the model-free extension later on.  


\section{A Primal Gradient Descent (GD) Approach}\label{sec:DP}

To deal with constraints in \eqref{eq:opt}, consider its  Lagrangian function by introducing the multiplier vector  $\lambda = \{\lambda_i\geq 0\}$, as
\begin{align}
	\mathcal{L}(K, \lambda) &=R_0(K)+ \textstyle \sum_{i \in \mathcal{I}} \lambda_i  [R_i(K)-c_i]\nonumber\\
	&=\lim_{T \rightarrow \infty} \frac{1}{T}\mathbb{E}\sum_{t=0}^{T-1} [x^\top_t Q_{\lambda} x_t + u^\top_t R_{\lambda} u_t]-c_{\lambda}
\end{align}
where we define  $Q_{\lambda} := Q+\sum_{i \in \mathcal{I}} \lambda_i Q_i $, and likewise for  $R_{\lambda}$ and $c_\lambda$. Clearly,  $\mathcal{L}(K, \lambda)$ shares the same structure as an unconstrained LQR cost which is suitable for first-order algorithms. {\KB For simplicity, consider that the problem \eqref{eq:opt} is feasible and thus $\lambda$ is finite \cite[Sec.~5.2]{boyd}. We consider the bounded set $\ccalY :=  [0, ~\Lambda]^{|\mathcal{I}|}$ for $\lambda$ with a large enough $\Lambda\in\mathbb{R}$, which can be set based on a feasible $K_0$.} Using the dual function $\ccalD(\lambda) : = \min_{K \in \mathcal{K}} \mathcal{L} (K, \lambda)$, the  dual problem becomes 
\begin{align}
	\max_{\lambda \in \mathcal{Y}} \ccalD(\lambda) = \max_{\lambda \in \mathcal{Y}}  \min_{K \in \mathcal{K}} \mathcal{L} (K, \lambda). \label{eq:DP}
\end{align}

As $\mathcal{L}(K, \lambda)$ is related to LQR cost, the inner minimization problem is not convex. Recent works \cite{malik,bu,harvard} have extensively analyzed the LQR cost which can be used to establish the local Lipschitz and smoothness properties of $\ccalL(K,\lambda)$. {Specifically, it is possible to find related constants that hold within a subset $\ccalG^0 \subset \ccalK$.} This compact sublevel set will be defined later on, but is first introduced here for bounding the constants as stated below.

\begin{lemma}[Lipschitz and smoothness]
	For any $\lambda$ and $K\in\ccalG^0$, the function $\ccalL(K,\lambda)$ is locally $L_0$-Lipschitz within a radius $\psi_K$; i.e., for $\forall K'\in \ccalG^0$ such that $\Vert K-K'\Vert \leq \psi_K$, we have $\Vert \ccalL(K,\lambda)-\ccalL(K',\lambda)\Vert \leq L_0 \Vert K-K'\Vert$. In addition, it is also locally $\ell_0$-smooth within a radius $\beta_K$, such that for $\forall K'\in \ccalG^0$ that satisfies $\Vert K-K'\Vert \leq \beta_K$, we have $\Vert \nabla\ccalL_{\ccalK}(K,\lambda)-\nabla\ccalL_{\ccalK}(K',\lambda)\Vert \leq \ell_0 \Vert K-K'\Vert$.
	 \label{lem:lips}
\end{lemma}

{Strictly speaking, the recent LQR analysis \cite{malik,harvard} asserts that Lipschitz and smoothness are only local properties, and thus the corresponding constants $L_K$ and $\ell_K$ depend on $K$. Nonetheless, using a compact set $\ccalG^0$,}
we can obtain the bounds {\KB that can hold for any $K \in \ccalG^0$, as given by}
{\KB
\begin{align}
	L_0:=\sup\limits_{K\in\mathcal{G}^0}{L_K},  ~\textrm{and}~\ell_0:=\sup\limits_{K\in\mathcal{G}^0}\ell_K. \label{eq:GDconst}
\end{align}
}
We can also determine a general neighborhood radius as
\begin{align}
   \rho_0:=\inf\limits_{K \in \ccalG^0} \min\{ \beta_K, \psi_K\} \label{eq:rho}
\end{align}
that holds for any $K\in\ccalG^0$ as well. 

Interestingly, the KKT conditions for problem \eqref{eq:DP} is related to the stationary point (SP) of a reformulated minimax problem. Recent results have shown that nonconvex-concave minimax problems can be solved using the so-termed Gradient Descent with max-oracle (GDmax) algorithm \cite{jin}. To this end, consider the problem
\begin{align}
 \min_{K\in\ccalK}   \Phi(K) 
~~\textrm{where} ~~\Phi(K):= \max_{\lambda \in \mathcal{Y}} \ccalL(K, \lambda) \label{eq:phi},
\end{align}
which is essentially the minimax  counterpart of problem \eqref{eq:DP}. As the Lagrangian function is linear in $\lambda$, it is possible to directly find the best $\lambda$ in \eqref{eq:phi}. Specifically, its $i$-th element, namely $\lambda_i$, depends on the feasibility of constraint $i$ under given $K$; i.e., $\lambda_i$ equals to $0$ if constraint $i$ is satisfied and  $\Lambda$ otherwise. Unfortunately, the function $\Phi(K)$ is not differentiable everywhere. To tackle this issue, we consider its \textit{Moreau envelope} $\Phi_\mu(\cdot)$ for a given $\mu > 0$, defined as
\begin{align}
    \Phi_\mu(K) := \min_{K'\in\ccalK} \Phi(K') + \frac{1}{2\mu}\Vert K'-K\Vert^2, ~~\forall K \in \ccalK. \label{eq:mor} 
\end{align}
It can be used for defining the SP of the non-differentiable $\Phi(K)$, following from \cite[Lemma 3.6]{GDA}.

\begin{lemma}
As $\ccalL(K,\lambda)$ is concave in $\lambda$ and $\ccalY$ is convex and bounded, Lemma \ref{lem:lips} asserts that $\Phi(K)$ is $\ell_0$-weakly convex and $L_0$-Lipschitz within the compact set $\ccalG^0$. Accordingly, its Moreau envelope $\Phi_{\mu_0}(K)$ is convex by setting $\mu_0 := 1/(2\ell_0)$. Hence, the $\epsilon$-SP of $\Phi(K)$, namely  $K_\epsilon$, satisfies  $\Vert \nabla\Phi_{\mu_0} (K_\epsilon) \Vert \leq \epsilon$. \label{lem:ep}
\end{lemma}

The properties of $\Phi(K)$  in Lemma \ref{lem:ep} follow from its relation to $\ccalL(K,\lambda)$, as detailed in \cite{GDA}. Even though it is non-differentiable, one can define the SP here based on $\Phi_{\mu_0}(K)$ which will be used for the convergence analysis of GD updates later on. Notably, the $\epsilon$-SP of $\Phi(K)$ is equivalently related to the  stationarity conditions for $\ccalL(K,\lambda)$. According to \cite[Prop.~4.12]{GDA}, one can utilize $K_\epsilon$ from Lemma~\ref{lem:ep} to generate the following pair $(\tilde{K}_\epsilon, \tilde{\lambda}_\epsilon)$  by performing an additional $O(\epsilon^{-2})$ number of gradient updates:
    \begin{align*}
        &\| \nabla_\ccalK \ccalL(\tilde{K}_\epsilon,\tilde{\lambda}_\epsilon) \| \leq \epsilon \\
        &\left\| \mathbb P_\ccalY\left(\tilde{\lambda}_\epsilon + (1/\ell_0)\nabla_\lambda \ccalL(\tilde{K}_\epsilon,\tilde{\lambda}_\epsilon)\right)-\tdlambda_\epsilon \right\| \leq \epsilon/\ell_0
    \end{align*}
where $\mathbb P_\ccalY$ stands for the projection onto $\ccalY$. 
{Clearly, when $\epsilon \rightarrow 0$ this  represents the Lagrangian optimality conditions for problem \eqref{eq:DP}, and thus the pair $(\tilde{K}_\epsilon,\tilde{\lambda}_\epsilon)$ can be viewed as the $\epsilon$-SP for $\ccalL(K,\lambda)$.}

We can solve \eqref{eq:phi} using iterative GD updates, as tabulated in Algorithm \ref{alg:GD}. With an initial $K^0$,  we need to find the subgradient of $\Phi(K^j)$ at every iteration $j$. Interestingly, this is equivalent to the gradient of $\ccalL$ over $K^j$ \cite{GDA}; i.e.,  
$\partial\Phi(K^{j})  = \nabla_\ccalK \mathcal{L}(K^{j}, \lambda^j) $
with $\lambda^j$ being the optimal multiplier for the given $K^{j}$. Hence, the Lagrangian $\ccalL$ will be used to perform the GD updates for $\Phi(K)$ minimization. The convergence of Algorithm \ref{alg:GD} can be established below,  with the detailed proof in Appendix~\ref{sec:APA}.

\begin{algorithm}[t]
	\SetAlgoLined
	\caption{Gradient Descent with max-oracle (GDmax)}
	\label{alg:GD}
	\DontPrintSemicolon
	{\bf Inputs:} A feasible policy $K^0$, upper bound $\Lambda$ for $\lambda$, threshold $\epsilon$, and the initial iteration index $j=0$.\; 
	Determine $L_0, \ell_0$, and $\rho_0$ using the set $\ccalG^0$ and compute the stepsize as in \eqref{eq:eta}.\;
	\While{$\|\nabla_{\ccalK}\ccalL(K^j, \lambda^j)\| > \epsilon $}
	{
		Obtain $\lambda^j \leftarrow \argmax_{\lambda \in \ccalY} \ccalL(K^j, \lambda)$\; 
		Update $K^{j+1} \leftarrow K^{j}-\eta \nabla_\ccalK \ccalL(K^j, \lambda^j)$;\;
		Set $j\leftarrow j+1$.
		}
	{\bf Return:} the final iterate $K^j$.
\end{algorithm}


\begin{theorem}
	\label{thm:GD}
With an initial $K^0 \in \ccalK$ and by setting stepsize  
\begin{align}
    \eta\leq\min \left \{\frac{\epsilon^2}{4\ell_0 L_0^2},\rho_0\right\}, \label{eq:eta}
\end{align}
Algorithm~\ref{alg:GD} is guaranteed to converge to $K_\epsilon$ for $\Phi(K)$, which can be used to obtain an $\epsilon$-SP for the dual problem \eqref{eq:DP}. {\KB The number of iterations required for attaining $K_\epsilon$ is  $O(\ell_0 L_2^2 \Phi_{\mu_0}(K^0)/\epsilon^4)$.} 
\end{theorem}

As discussed in Appendix~\ref{sec:APA}, we can  bound the iterative changes in $\Phi_{\mu_0}(K^j)$, which ensures that the sequence $\{\Phi_{\mu_0}(K^j)\}$ is  non-increasing. Thus, if we define the sublevel set to be
\begin{align}
\ccalG^0 :=\{K \in \ccalK \vert \Phi_{\mu_0}(K) \leq \Phi_{\mu_0}(K^0) \}, \label{eq:G0}
\end{align}
then the iterates $\{K^j\}$ are guaranteed to be within  $\ccalG^0$. This is exactly how one can bound the constants $L_0$ and $\ell_0$ as given by \eqref{eq:GDconst}. Of course,  the choice of $\mu_0$ in the sublevel set $\ccalG_0$ depends on $\ell_0$, which may not be known before $\ccalG^0$ is constructed. This issue is discussed in the following remark.  


\begin{remark}[Sublevel set] \label{re:sub}
With initial $K^0$ given, the set $\ccalG^0$ is defined with the value $\mu_0$, which depends on  the upper bound of $\ell_K$ within $\ccalG^0$ as shown in \eqref{eq:GDconst}. This dependence can be addressed by determining the value of $\mu_0$ in an adaptive fashion. Starting with a rough estimate of $\ell_0$ and $\mu_0$, one can first construct a $\ccalG^0$ and compare the resultant bound with the original estimate on $\ell_0$. If the latter is larger, then $\ccalG^0$ works well. Otherwise, one can gradually increase the $\ell_0$ estimate to achieve that condition. Our experimental experience suggests some conservative choice of stepsize 
can ensure the convergence in practice. 
\end{remark}



\section{Stochastic GD for Model-free Learning} \label{sec:SG}

To account for unknown system dynamics, we extend the GDmax approach to a model-free setting. The iterative gradient will be obtained via the zero-order optimization \cite{ZOO}. Unfortunately, this stochastic gradient update can complicate the convergence analysis as detailed later, mainly due to the aforementioned issue on local  properties of LQR cost.

Zero-order policy gradient (ZOPG) has been popularly developed in recent years for model-free gradient-based learning. It provides an unbiased gradient estimate in an efficient manner. For the function $\Phi(K)$, ZOPG aims to evaluate the function value at any $K$ under a structured, random perturbation from the set $\mathcal{S}_{\ccalK} = \{U \in \ccalK: \Vert U \Vert = 1\}$, as detailed in Algorithm~\ref{alg:ZOO}. Note that the structure of perturbation $U$ is the same to that of $K$ with non-zero entries randomly sampled from e.g., the uniform distribution, followed by a  normalization step to ensure unity norm. Given a smoothing radius $r>0$, the ZOPG is estimated using the resultant $\Phi(K+rU)$ from this perturbation by finding the corresponding optimal $\lambda$ in \eqref{eq:phi}. We denote $n_\ccalK$ as the total number of nonzero entries in $\ccalK$, which is used to scale the gradient estimate. Since the estimated $\hat{\nabla}_{\ccalK}\ccalL$ follows from matrix $U$, it maintains the same sparse structure given by $\ccalK$. 

\begin{algorithm}[t]
	\SetAlgoLined
	\caption{Zero-Order Policy Gradient (ZOPG)}
	\label{alg:ZOO}
	\DontPrintSemicolon
	{\bf Inputs:} smoothing radius $r$, 
	the policy $K$ and its perturbation  $U \in \mathcal{S}_{\mathcal{K}}$, both of $n_\ccalK$ non-zeros.\; 
	Obtain $\lambda' \leftarrow \argmax_{\lambda \in \ccalY} \ccalL(K+rU, \lambda)$;\; 
Estimate the gradient $\hat{\nabla}_{\ccalK} \ccalL (K; U) = \frac{n_\ccalK}{r} \ccalL(K+rU, \lambda') U$.\;
	{\bf Return:} $\hat{\nabla}_{\ccalK} \ccalL (K; U)$.
\end{algorithm}


{The stochastic ZOPG will make it more difficult to maintain  the  iterative updates to stay within a sublevel set, and likewise for bounding Lipschitz and smoothness constants. 
Fortunately, \cite{malik} has developed an approach to attain this condition with a high probability. Specifically, one can set up a ten-fold sublevel set, given by  
\begin{align}
    \ccalG^1 := \{K\in\ccalK |\Phi_{\mu_0}(K)\leq 10~\Phi_{\mu_0}(K^0)  \}. \label{eq:ten}
\end{align}
{\KB Using $\ccalG^1$, one can determine $L_0, \ell_0,$ and $\rho_0$ over the set $\ccalG^1$ similar to \eqref{eq:GDconst}-\eqref{eq:rho}, and they will be used for the convergence analysis.} 
Note that the choice of $\mu_0$ in $\ccalG^1$ depends on the $\ell_0$ value,  which can be addressed as discussed in Remark \ref{re:sub}.

Algorithm~\ref{alg:SGD} tabulates the ZOPG-based model-free learning approach for solving \eqref{eq:DP}, termed as the Stochastic Gradient Descent with max-oracle (SGDmax) \cite{jin}. Its convergence guarantee can be established with the detailed proof in Appendix~\ref{sec:APB}. 


\begin{algorithm}[t]
	\SetAlgoLined
	\caption{Stochastic Gradient Descent with max-oracle (SGDmax)}
	\label{alg:SGD}
	\DontPrintSemicolon
	{\bf Inputs:} A feasible policy $K^0$, upper bound $\Lambda$ for $\lambda$, threshold $\epsilon$, and number of ZOPG samples $M$.\; 
		Determine $L_0, \ell_0$, and $\rho_0$ with the set $\ccalG^1$ and compute $r$, $\eta$, and $J$ as in \eqref{eq:thm2}; \;
	\For{$j = 0, 1, \ldots, J-1$}{
		\For{$s=1,\ldots,M$}{
		Sample the random $U_s \in \ccalS_\ccalK$; \;
		Use Algorithm~\ref{alg:ZOO} to return $\hat{\nabla}\ccalL_\ccalK(K^j;U_s)$. 
		}
		Update $K^{j+1}\!\leftarrow\! K^{j}\!-\!\eta \left(\frac{1}{M} \sum_{s=1}^M \hat{\nabla} \ccalL(K^j;U_s)\right)$.}
	{\bf Return:} the final iterate $K^{J}$.
\end{algorithm}

\begin{theorem}
	\label{thm:SGD} 
 With an initial $K^0 \in \ccalK$  and a given $\epsilon > 0$, we can set the parameters as 
	\begin{align}
		&r \leq \min\Big\{\rho_0, \frac{L_0 \sqrt{M}}{\ell_0} \Big\}, \; \eta \leq \frac{\epsilon^2}{\alpha\ell_0(L_0^2+\ell_0^2r^2/M)},\nonumber\\
		\textrm{and}~&J = \frac{2\sqrt{10\alpha}\Phi_{\mu_0}(K^0)}{\eta \epsilon^2}\label{eq:thm2}
	\end{align}
with $L_0, \ell_0$ and $\rho_0$ being specified using $\ccalG^1$, and {a large constant $\alpha$.} This way, Algorithm \ref{alg:GD} converges to the $\epsilon$-SP $K_\epsilon$ with probability of at least  $(0.9-\frac{4}{\alpha}-\frac{4}{\sqrt{10\alpha}})$.
\end{theorem}
}

Last, the proposed algorithms can be easily extended to the case of full feedback $K$, with computational advantages over existing solutions as discussed below.

\begin{remark}[Full feedback $K$] \label{re:full}
For the full feedback case, we can directly implement the proposed Algorithms \ref{alg:GD}-\ref{alg:SGD} by dropping the structured set $\ccalK$. This setting has been considered in \cite{zhao} by using a dual-ascent based double-loop scheme where the inner-loop minimizes $K$ till convergence for any fixed $\lambda$.  In contrast, our proposed algorithms eliminate this inner-loop, which is more computationally efficient. Investigating the global convergence property of our proposed SGDmax algorithm for the full feedback case constitutes as an interesting future direction. 
\end{remark}

\section{Numerical Tests}\label{sec:NT}
\begin{figure}[t]
	\centering
	\includegraphics[width=0.5\linewidth]{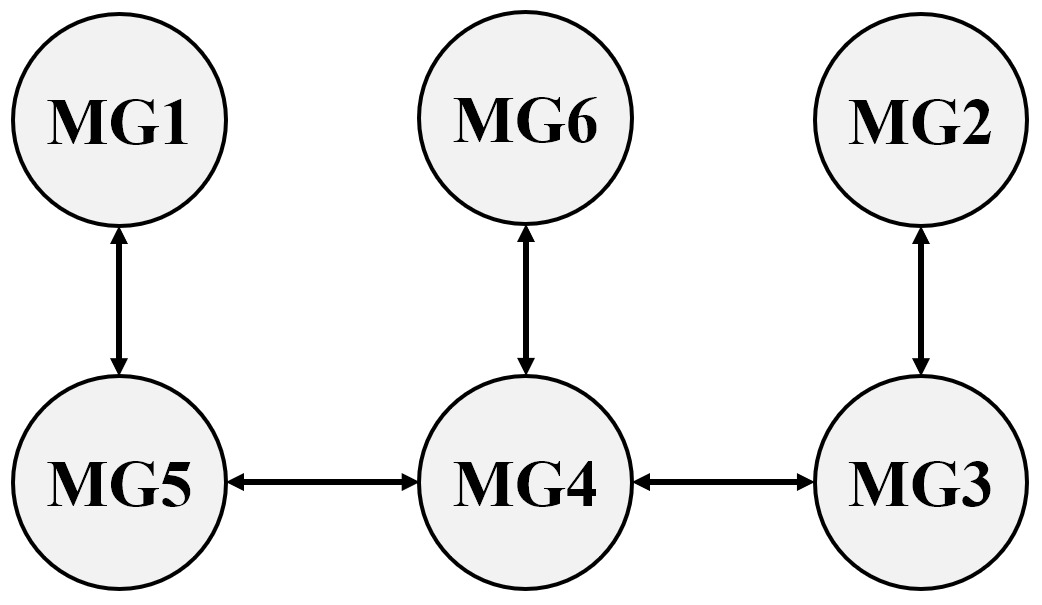}
	\caption{A radially connected networked microgrid system.}
	\label{fig:network}
\end{figure}

To demonstrate the effectiveness of the proposed model-free learning approach, we consider the load frequency control (LFC) problem in a low-inertia networked microgrid (MG) system with a risk constraint on the frequency states.  Fig.~\ref{fig:network} depicts a radially connected system with $N=6$ MGs, while Table~\ref{tb:parameter} lists the model information which follows from \cite{DLQR}. Consider the communication graph to be the same as the MG network show in Fig.~\ref{fig:network}. Thus, each MG~$a$ can only exchange information with their neighboring MGs that are physically connected by tie-lines, and the structured feedback $\ccalK$ is specified accordingly.

\begin{table}[t]
	\centering
	\caption{List of parameter and their values}
	\label{tb:parameter}
	\begin{tabular}{c c c c}
		\toprule
		\textbf{Parameter} & \text{Symbol} & \text{Value} & \text{Units}\vspace*{3pt} \\ \hline \hline
		Damping Factor & $D$ & 16.66 & MW/Hz\\
		Speed Droop & $R$ & $1.2 \times 10^{-3}$ & Hz/MW\\
		Turbine Static Gain & $K_{t}$ & 1 & MW/MW\\
		Turbine Time Constant & $T_{t}$ & 0.3 & s\\
		Area Static Gain & $K_{p}$ & 0.06 & Hz/MW\\
		Area Time Constant & $T_{p}$ & 24 & s\\
		Tie-line Coefficient & $K_{tie}$ & 1090 & MW/Hz\\ \bottomrule
	\end{tabular}
\end{table}


Each MG~$a$ is assumed to follow linearized power-frequency dynamics including turbine swing and primary control based on the automatic generation control (AGC) signal. Thus, the following symbols all correspond to the deviation from steady-state values as denoted by $\Delta$, with the parameters listed in Table~\ref{tb:parameter}. First, the primary frequency control in each MG~$a$ is proportional to frequency deviation as $\Delta P_{f,a} = -(1/R_a) \Delta f_a$ based on the given droop $R_a$. Second, the secondary AGC signal $\Delta P_{C,a}$ constitutes as the control action $u_t$ in \eqref{eq:dynamics} to be designed. The two controls jointly determine the power output of MG $a$ as denoted by $\Delta P_{G,a}$. Last, $\Delta f_a$ is also affected by the unknown load demand deviation $\Delta P_{L,a}$ and the total power inflow $\Delta P_{tie,a}$, in addition to $\Delta P_{G,a}$. Note that $\Delta P_{tie,a}$ is the total tie-line power inflow from all neighboring MGs due to their frequency differences, as 
\begin{align}
	\Delta P_{tie,a} = \int \sum_{a \leftrightarrow b} K_{tie,a} (\Delta f_a - \Delta f_b) dt,
\end{align}   
where $a \leftrightarrow b$ indicates two MGs are connected to each other. In addition to the MG dynamics, the Area Control Error (ACE) defined as $z_a := \beta_a \Delta f_a + \Delta P_{tie,a}$ is also a state variable as an integral control input with the bias factor $\beta_a=D_a + 1/R_a$  \cite{DLQR_en}.   



Hence, MG $a$ has the state vector $x_a = [\Delta f_a, \Delta P_{G,a},$ $\Delta P_{tie,a},\int z_a]^\top$ and the control action  $u_a = \Delta P_{C,a}$, with load disturbance $w_a = \Delta P_{L,a}$. Assuming all MGs having the same parameter values, we can drop the parameter index $a$ and represent the aggregated network dynamics by:
\begin{align*}
	\dot{{x}} = (I_N \otimes A_1 +L \otimes A_2){x} + (I_N \otimes B_u) {u} + (I_N \otimes B_w)\tilde{w} 
\end{align*}
with each variable collecting all MGs' respective state, action, and disturbance. In addition, the system matrices are given by
%
%
\begin{align*}
	&A_1 = \begin{bmatrix}
		-\frac{1}{T_{p}} & \frac{K_{p}}{T_{p}} & -\frac{K_{p}}{T_{p}} & 0\\
		-\frac{K_{t}}{R T_{t}} & -\frac{1}{T_{t}} & 0 & 0\\
		0 & 0 & 0 & 0\\
		\beta & 0 & 1 & 0
	\end{bmatrix}, \nonumber \\ 
	&A_2 = \begin{bmatrix}
		0 & 0 & 0 & 0\\
		0 & 0 & 0 & 0\\
		K_{tie} & 0 & 0 & 0\\
		0 & 0 & 0 & 0
	\end{bmatrix},
	B_u = \begin{bmatrix}
		0\\		\frac{K_t}{T_t}	\\0	\\0
	\end{bmatrix},
	B_w = \begin{bmatrix}
		-\frac{K_p}{T_p}\\	0\\	0\\	0
	\end{bmatrix}
\end{align*}

For the aggregated dynamics, the LQR objective cost is specified by 
\begin{align*}
{Q} = I_{N_L} \otimes Q_a,~\textrm{and}~{R} = I_{N_L} \otimes R_a
\end{align*} 
where the matrices $Q_a$ and $R_a$ are same for every MG $a$ and aim to penalize the deviation of both state and action from steady-state values. As discussed in Section~\ref{sec:PF}, we further consider a risk constraint $R_c(\cdot)$ in \eqref{eq:objr2} for reducing the mean-variance risk in order to improve frequency regulation. 
%

\begin{figure}[t]
	\centering
	\includegraphics[width=0.8\linewidth]{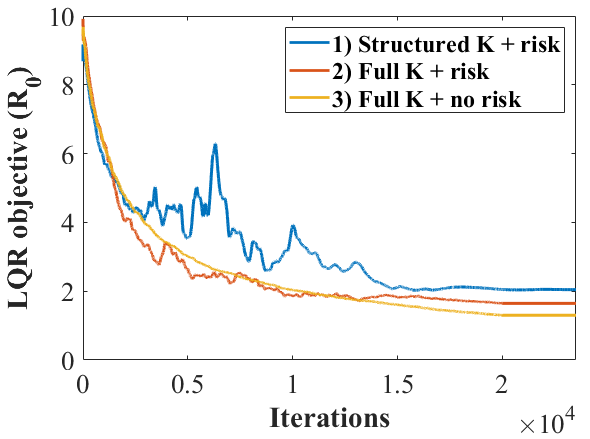}
	\caption{Comparison of LQR objective trajectories for the three cases.}
	\captionsetup{justification=centering}
	\label{fig:conv}
\end{figure}

We consider the following three cases to demonstrate the impact of structured $K$ along with the risk constraint: 
\begin{itemize}
	\item Case 1): Structured $K$ with risk constraint
	\item Case 2): Full $K$ with risk constraint
	\item Case 3): Full $K$ without risk constraint
\end{itemize}
For cases 1 and 2, we implemented Algorithm \ref{alg:SGD} using SGDmax while a simple ZOPG-based algorithm \cite{malik} was used for case 3. {For all algorithms, we picked a small stepsize of $\eta = 10^{-4}$ with a smoothing radius $r=1$ and $M=100$ samples for ZOPG.} All three cases have shown to converge to a steady-state with sufficient updates, as shown by Fig.~\ref{fig:conv}. In particular, the LQR cost attained by case 2 is slightly over that by case 3, suggesting a global convergence result for SGDmax in full feedback case as discussed in Remark \ref{re:full}.  Case 1 demonstrates the highest steady-state LQR cost out of the three, as it has the most restrictive conditions.  However, the minimum LQR cost by case 1 is still pretty close to that by case 3, implying some good optimality gap. Notably, case 1 has shown some large fluctuations along the learning process, indicating a complicated geometry that the problem may have.

\begin{figure}
	\centering
	\subfigure[]{\label{fig:freq1}\includegraphics[width=\linewidth]{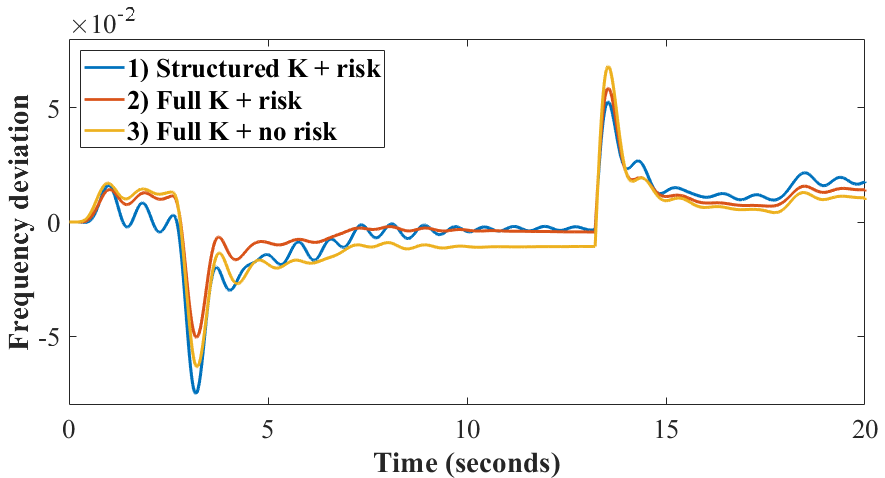}}
	\subfigure[]{\label{fig:tie1}\includegraphics[width=\linewidth]{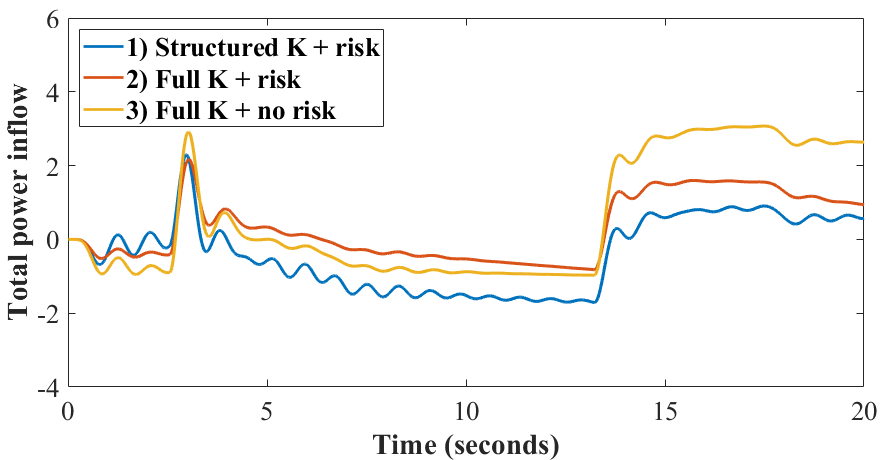}}
	\caption{Comparison of the (a) frequency deviation and (b) total power inflow at MG~$2$ for the three cases.}
	\captionsetup{justification=centering}
	\label{fig:total}
\end{figure}

We also test the converged policy by each case by generating a scenario that all six MGs have some random load changes in a 20-second window. Each area experiences a step load change at a random time. Fig.~\ref{fig:total} compares the frequency deviation and the total power inflow for MG~$2$. Clearly, Fig.~\ref{fig:total}(a) demonstrates that the risk constraint can effectively reduce the frequency deviation, as case 2 has the smallest deviation among all three. With the risk constraint, case 1 tends to exhibit great frequency performance as well, but also shows some small oscillations possibly due to the structured feedback policy. This observation points out that limited information exchange can potentially affect the control performance. Similar patterns have been observed in Fig.~\ref{fig:total}(b). While case 1 can maintain the tie-line inflow at the same level as case 2, it still has more noticeable oscillations. As the power inflow is proportional to frequency difference, reducing the risk of frequency deviation can enhance the performance in maintaining the level of power inflow. 



To sum up, our numerical tests have validated the convergence performance of the proposed SGDmax based policy gradient method for risk-constrained LQR problem with structured policy. The effectiveness of risk constraint in mitigating large state deviation have been verified, while the sparse structure of $K$ has shown to save communication overhead at the cost of  transient oscillations.


\section{Conclusions}\label{sec:CN}
The paper developed a model-free learning framework for risk-constrained LQR problem under structured feedback in a networked setting. By dualizing the risk constraint, we consider the minimax reformulation of the dual problem and leverage the stochastic (S)GDmax algorithms to approach the stationary points (SPs). Specifically, the SGDmax algorithm relies on the ZOPG-based updates, making it suitable for model-free learning. Using the recent results on the local Lipschitz and smoothness of LQR cost, convergence of the (S)GDmax algorithms can be established by properly bounding the related constants for choosing the stepsize. Notably, for SGDmax the convergence can only be shown with a high probability, due to the additional noise in the gradient estimate. Numerical tests on a networked microgrid system  have validated the convergence of our proposed algorithms while demonstrating the impact of risk and structured constraints for the LQR problem. Exciting future research directions open up on investigating the landscape for the converged SP in the structured feedback case and establishing the global convergence for the full feedback case.
%
\bibliography{bibliography}

\begin{thebibliography}{10}
\providecommand{\url}[1]{#1}
\csname url@samestyle\endcsname
\providecommand{\newblock}{\relax}
\providecommand{\bibinfo}[2]{#2}
\providecommand{\BIBentrySTDinterwordspacing}{\spaceskip=0pt\relax}
\providecommand{\BIBentryALTinterwordstretchfactor}{4}
\providecommand{\BIBentryALTinterwordspacing}{\spaceskip=\fontdimen2\font plus
\BIBentryALTinterwordstretchfactor\fontdimen3\font minus
  \fontdimen4\font\relax}
\providecommand{\BIBforeignlanguage}[2]{{%
\expandafter\ifx\csname l@#1\endcsname\relax
\typeout{** WARNING: IEEEtran.bst: No hyphenation pattern has been}%
\typeout{** loaded for the language `#1'. Using the pattern for}%
\typeout{** the default language instead.}%
\else
\language=\csname l@#1\endcsname
\fi
#2}}
\providecommand{\BIBdecl}{\relax}
\BIBdecl

\bibitem{LQR1}
R.~E. Kalman \emph{et~al.}, ``{Contributions to the theory of optimal
  control},'' \emph{Bol. soc. mat. mexicana}, vol.~5, no.~2, pp. 102--119,
  1960.

\bibitem{LQR2}
B.~D.~O. Anderson and J.~B. Moore, \emph{{Optimal Control: Linear Quadratic
  Methods}}.\hskip 1em plus 0.5em minus 0.4em\relax USA: Prentice-Hall, Inc.,
  1990.

\bibitem{bu}
J.~Bu, A.~Mesbahi, M.~Fazel, and M.~Mesbahi, ``{LQR through the lens of first
  order methods: Discrete-time case},'' \emph{arXiv preprint arXiv:1907.08921},
  2019.

\bibitem{malik}
D.~Malik, A.~Pananjady, K.~Bhatia, K.~Khamaru, P.~Bartlett, and M.~Wainwright,
  ``{Derivative-free methods for policy optimization: Guarantees for linear
  quadratic systems},'' in \emph{The 22nd International Conference on
  Artificial Intelligence and Statistics}.\hskip 1em plus 0.5em minus
  0.4em\relax PMLR, 2019, pp. 2916--2925.

\bibitem{tsiamis}
A.~Tsiamis, D.~S. Kalogerias, L.~F.~O. Chamon, A.~Ribeiro, and G.~J. Pappas,
  ``{Risk-Constrained Linear-Quadratic Regulators},'' in \emph{{2020 59th IEEE
  Conference on Decision and Control (CDC)}}, 2020, pp. 3040--3047.

\bibitem{tsiamis2}
A.~Tsiamis, D.~S. Kalogerias, A.~Ribeiro, and G.~J. Pappas, ``{Linear Quadratic
  Control with Risk Constraints},'' \emph{arXiv preprint arXiv:2112.07564},
  2021.

\bibitem{zhao}
F.~Zhao, K.~You, and T.~Ba{\c{s}}ar, ``{Global Convergence of Policy Gradient
  Primal-dual Methods for Risk-constrained LQRs},'' \emph{arXiv preprint
  arXiv:2104.04901}, 2021.

\bibitem{paternain2019constrained}
S.~Paternain, L.~Chamon, M.~Calvo-Fullana, and A.~Ribeiro, ``Constrained
  reinforcement learning has zero duality gap,'' \emph{Advances in Neural
  Information Processing Systems}, vol.~32, 2019.

\bibitem{ding2020natural}
D.~Ding, K.~Zhang, T.~Basar, and M.~Jovanovic, ``Natural policy gradient
  primal-dual method for constrained markov decision processes,''
  \emph{Advances in Neural Information Processing Systems}, vol.~33, pp.
  8378--8390, 2020.

\bibitem{shah2013cal}
P.~Shah and P.~A. Parrilo, ``H2-optimal decentralized control over posets: A
  state-space solution for state-feedback,'' \emph{IEEE Transactions on
  Automatic Control}, vol.~58, no.~12, pp. 3084--3096, 2013.

\bibitem{DEC}
M.~Rotkowitz and S.~Lall, ``{A Characterization of Convex Problems in
  Decentralized Control},'' \emph{IEEE Transactions on Automatic Control},
  vol.~51, no.~2, pp. 274--286, 2006.

\bibitem{ye2021sample}
L.~Ye, H.~Zhu, and V.~Gupta, ``On the sample complexity of decentralized linear
  quadratic regulator with partially nested information structure,''
  \emph{arXiv preprint arXiv:2110.07112}, 2021.

\bibitem{dorfler2013sparse}
F.~D{\"o}rfler, M.~R. Jovanovi{\'c}, M.~Chertkov, and F.~Bullo, ``Sparse and
  optimal wide-area damping control in power networks,'' in \emph{2013 American
  Control Conference}.\hskip 1em plus 0.5em minus 0.4em\relax IEEE, 2013, pp.
  4289--4294.

\bibitem{chakrabortty2021wide}
A.~Chakrabortty, ``Wide-area control of power systems: Employing data-driven,
  hierarchical reinforcement learning,'' \emph{IEEE Electrification Magazine},
  vol.~9, no.~1, pp. 45--52, 2021.

\bibitem{DCM2}
C.~Wang, J.~Duan, B.~Fan, Q.~Yang, and W.~Liu, ``{Decentralized
  High-Performance Control of DC Microgrids},'' \emph{IEEE Transactions on
  Smart Grid}, vol.~10, no.~3, pp. 3355--3363, 2019.

\bibitem{DLQR}
E.~E. Vlahakis, L.~D. Dritsas, and G.~D. Halikias, ``{Distributed LQR design
  for identical dynamically coupled systems: Application to Load Frequency
  Control of multi-area Power Grid},'' in \emph{2019 IEEE 58th Conference on
  Decision and Control (CDC)}.\hskip 1em plus 0.5em minus 0.4em\relax IEEE,
  2019, pp. 4471--4476.

\bibitem{harvard}
Y.~Li, Y.~Tang, R.~Zhang, and N.~Li, ``{Distributed reinforcement learning for
  decentralized linear quadratic control: A derivative-free policy optimization
  approach},'' \emph{IEEE Transactions on Automatic Control}, 2021.

\bibitem{lavaeri}
H.~Feng and J.~Lavaei, ``{On the exponential number of connected components for
  the feasible set of optimal decentralized control problems},'' in \emph{2019
  American Control Conference (ACC)}.\hskip 1em plus 0.5em minus 0.4em\relax
  IEEE, 2019, pp. 1430--1437.

\bibitem{GDA}
T.~Lin, C.~Jin, and M.~Jordan, ``{On gradient descent ascent for
  nonconvex-concave minimax problems},'' in \emph{International Conference on
  Machine Learning}.\hskip 1em plus 0.5em minus 0.4em\relax PMLR, 2020, pp.
  6083--6093.

\bibitem{ZOO}
J.~C. Spall, ``{A one-measurement form of simultaneous perturbation stochastic
  approximation},'' \emph{Automatica}, vol.~33, no.~1, pp. 109--112, 1997.

\bibitem{boyd}
S.~Boyd and L.~Vandenberghe, \emph{{Convex optimization}}.\hskip 1em plus 0.5em
  minus 0.4em\relax Cambridge university press, 2004.

\bibitem{jin}
C.~Jin, P.~Netrapalli, and M.~Jordan, ``{What is local optimality in
  nonconvex-nonconcave minimax optimization?}'' in \emph{International
  Conference on Machine Learning}.\hskip 1em plus 0.5em minus 0.4em\relax PMLR,
  2020, pp. 4880--4889.

\bibitem{DLQR_en}
E.~Vlahakis, L.~Dritsas, and G.~Halikias, ``{Distributed LQR design for a class
  of large-scale multi-area power systems},'' \emph{Energies}, vol.~12, no.~14,
  p. 2664, 2019.

\end{thebibliography}
\bibliographystyle{IEEEtran}
\itemsep2pt

%
\numberwithin{equation}{subsection}
\numberwithin{definition}{subsection}
\numberwithin{lemma}{subsection}

\appendix
\subsection{Proof of Theorem~\ref{thm:GD}} \label{sec:APA}

The key step is to ensure that the iterates stay within the sublevel set $\ccalG^0$ defined in Section \ref{sec:DP}. To this end, consider the function $\Phi(\cdot)$ in \eqref{eq:phi} with its \textit{Moreau envelope} $\Phi_{\mu}(\cdot)$ defined as \eqref{eq:mor}. 
Based on Lemma~\ref{lem:ep}, the problem becomes to show the convergence of $\Phi_{\mu_0}(\cdot)$ instead. To bound the iterative change in $\Phi_{\mu_0}(\cdot)$, one can use $L_0$-Lipschitz and $\ell_0$-weakly convex properties of $\Phi(\cdot)$ to analyze the update $K^{j+1} \leftarrow K^{j}-\eta \nabla\Phi(K^j)$ and obtain \cite[Lemma D.3]{GDA}
\begin{align}
	\Phi_{\mu_0}(K^{j+1}) &\leq \Phi_{\mu_0}(K^{j}) - \frac{\eta}{4}\left \Vert \nabla \Phi_{\mu_0}(K^{j}) \right \Vert^2 + \eta^2 \ell_0 L_0^2. \label{eq:pbd}
\end{align}
{Note that $\eta \leq \rho_0$ is needed to apply the constants $L_0$ and $\ell_0$. Furthermore, by setting $\eta\leq\epsilon^2/(4\ell_0 L_0^2)$, the last term is upper bounded by $\epsilon^4/(16\ell_0 L_0^2)$, while the second term is lower bounded by the same value as $\Vert \nabla \Phi_{\mu_0}(K^j) \Vert > \epsilon$ holds before reaching $K_\epsilon$. Therefore, we can guarantee that $\Phi_{\mu_0}(K^j)$ is non-increasing and $K^j \in \ccalG^0 \; \forall j$. As a result, $L_0$-Lipschitz and $\ell_0$-smoothness properties hold throughout the iterative updates.} 

To verify  the SP condition in Lemma \ref{lem:ep},  summing up \eqref{eq:pbd} over $j = 0, 1, \ldots, J-1$ yields
\begin{align*}
	\frac{1}{J}\sum_{j=0}^{J-1} \left\Vert \nabla\Phi_{\mu_0}(K^j)\right\Vert^2 &\leq \frac{4\!\left[\Phi_{\mu_0}(K^0)\!-\!\Phi_{\mu_0}(K^{J})\right]}{J\eta}\!+\!4\eta \ell_0 L_0^2\\
	&\leq \frac{4\ell_0 L_0^2\Phi_{\mu_0}(K^0)}{J\epsilon^2}\!+\!\epsilon^2
\end{align*}
where the second step uses the choice of stepsize in \eqref{eq:eta}. {As $K^0$ is stable, the value $\Phi_{\mu_0}(K^0)$ is finite and thus the first term  is in the order of $\epsilon^2$ with $J= O(\ell_0 L_2^2 \Phi_{\mu_0}(K^0) /\epsilon^4)$ iterations. As a result, the gradient norm $\Vert \nabla\Phi_{\mu_0}(K^j)\Vert$ eventually approaches $\epsilon$, satisfying the $\epsilon$-SP condition. } \hfill \qedsymbol

\subsection{Proof of Theorem~\ref{thm:SGD}}\label{sec:APB}

Similar to Appendix~\ref{sec:APA}, the key lies in the iterative analysis of function $\Phi_{\mu_0}(K)$, or in this case its expectation. First, due to the noisy gradient of ZOPG, one can obtain the following  inequality similar to \eqref{eq:pbd} \cite[Lemma D.4]{GDA}:
\begin{align}
	\mathbb{E}&\left[\Phi_{\mu_0}(K^{j+1})\right]\leq \mathbb{E}\left[\Phi_{\mu_0}(K^{j})\right]\nonumber\\
	& - \frac{\eta}{4} \mathbb{E} \Vert \nabla \Phi_{\mu_0}(K^{j})\Vert^2 + \eta^2 \ell_0 \left(L_0^2+{\ell_0^2 r^2}/{M}\right) \label{eq:pbd2}
\end{align}
where the last term is because the noise variance of each ZO gradient sample with a smoothing radius $r$ is bounded by $\ell_0^2r^2$ as shown in \cite{malik,harvard}, while $M$ is the total number of samples. {Note that by choosing the smoothing radius $r$ as in Theorem~\ref{thm:SGD}, we prevent the overall noise variance $(\ell_0^2r^2/M)$ to be dominant in the last term.} Moreover, $r$ needs to be smaller than $\rho_0$ to ensure that each ZOPG iteration can use the local Lipschitz and smoothness constants. 



Summing up \eqref{eq:pbd2} over iterations $j=0,\ldots,J-1$ yields
\begin{align}
    &\frac{1}{J}\sum_{j=0}^{J-1}\mathbb{E}\left[\Vert \nabla\Phi_{\mu_0}(K^j)\Vert^2\right] \nonumber\\
    \leq&\frac{4\left[\Phi_{\mu_0}(K^0)\!-\!\mathbb{E}[\Phi_{\mu_0}(K^J)]\right]}{J\eta}+\!4\eta\ell_0(L_0^2+\ell_0^2r^2/M). \nonumber
\end{align}
With $\eta=O(\epsilon^2)$ and $J$ inversely proportional to $\eta\epsilon^2$ given in Theorem~\ref{thm:SGD}, this upper bound is in the order of $\epsilon^2$, as detailed soon. To eliminate the expectation therein, one can analyze the probability of exceeding $\epsilon^2$ by  considering whether $\{K^j\}$ exceeds $\ccalG^1$ within $J$ iterations, as given  by
%
\begin{align}
    &\mathbb{P}\left( \frac{1}{J} \sum_{j=0}^{J-1} \Vert \nabla\Phi_{\mu_0}(K^j)\Vert^2 \geq \epsilon^2 \right) \nonumber\\
    = & \mathbb{P}\left( \frac{1}{J} \sum_{j=0}^{J-1}\Vert \nabla\Phi_{\mu_0}(K^j)\Vert^2 \geq \epsilon^2, \tau > J \right) \nonumber\\
    + &\mathbb{P}\left( \frac{1}{J} \sum_{j=0}^{J-1}\Vert \nabla\Phi_{\mu_0}(K^j)\Vert^2 \geq \epsilon^2, \tau \leq J \right), \label{eq:prob}
\end{align}
where $\tau := \min\{j \geq 0: K^j \notin \ccalG^1\}$. The first term of \eqref{eq:prob} can be bounded by
\begin{align}
&\mathbb{P}\left( \frac{1}{J} \sum_{j=0}^{J-1} \Vert \nabla\Phi_{\mu_0}(K^j)\Vert^2 \geq \epsilon^2, \tau > J \right) \nonumber\\
\leq & \frac{1}{\epsilon^2}\mathbb{E}\left[\frac{1}{J} \sum_{j=0}^{J-1} \Vert \nabla\Phi_{\mu_0}(K^j)\Vert^2\right] \nonumber\\
\leq & \frac{1}{\epsilon^2}\left\{\frac{4\left[\Phi_{\mu_0}(K^0)\!-\!\mathbb{E}[\Phi_{\mu_0}(K^J)]\right]}{J\eta} + 4\eta\ell_0(L_0^2+\ell_0^2r^2/M)\right\}\nonumber\\
\leq &\frac{4\Phi_{\mu_0}(K^0)}{J\eta\epsilon^2} + \frac{4}{\alpha} = \frac{4}{\beta} + \frac{4}{\alpha}  \label{eq:prob1}
\end{align}
where the first step follows from the Markov's inequality, while the last one uses the parameter settings in \eqref{eq:thm2} {with $\beta$ simplifying the first fractional term to be determined soon}. 

In addition, the second term can be bounded by recognizing  that the sequence $Y^j:= \Phi_{\mu_0}(K^{\min(j,\tau)})+(J-j)\eta\ell_0(L_0^2+\ell_0^2r^2/M)$ is a supermartingale, as shown in \cite{malik}. Thus, using the Doob's maximal inequality for supermartingales, one can bound the second term of \eqref{eq:prob} as
\begin{align}
    &\mathbb{P}\left( \frac{1}{J} \sum_{j=0}^{J-1} \Vert \nabla\Phi_{\mu_0}(K^j)\Vert^2 \geq \epsilon^2, \tau \leq J \right)\nonumber\\
    \leq &\mathbb{P}(\tau \leq J)\nonumber\\ 
    \leq &\frac{\Phi_{\mu_0}(K^0)+J\eta^2\ell_0(L_0^2+\ell_0^2r^2/M)}{10\Phi_{\mu_0}(K^0)}\nonumber\\
    \leq &\frac{1}{10} + \frac{J\eta\epsilon^2}{10\alpha\Phi_{\mu_0}(K^0)} = \frac{1}{10} +\frac{\beta}{10\alpha} \label{eq:prob2}
\end{align}
where the first step relaxes the probability, the second step follows from Doob's maximal inequality, while the last one again uses the parameter settings in \eqref{eq:thm2}. Therefore, the probability that $\{K^j\}$ exceeds $\ccalG^1$ before $J$ is bounded, and we can ensure that $\{K^j\}$ is within $\ccalG^1$ with a high probability. This is why the compact sublevel set $\ccalG^1$ can be used to bound the Lipschitz and smoothness constants for $\Phi(K)$.

By substituting \eqref{eq:prob1}-\eqref{eq:prob2}, the overall probability becomes
\begin{align*}    
    \mathbb{P}\left( \frac{1}{J} \sum_{j=0}^{J-1} \Vert \nabla\Phi_{\mu_0}(K^j)\Vert^2 \geq \epsilon^2 \right)&\leq \frac{1}{10}+ \frac{4}{\beta} + \frac{\beta}{10\alpha} + \frac{4}{\alpha}\\
    &\leq \frac{1}{10}+ \frac{4}{\alpha} + \frac{4}{\sqrt{10\alpha}}
\end{align*}
where the last step uses the best choice of $\beta = 2\sqrt{10\alpha}$. As a result, with probability $(0.9-\frac{4}{\alpha}-\frac{4}{\sqrt{10\alpha}})$, the $\epsilon$-SP can be attained by the iterations $\{K^j\}$ within $J$ iterations. Note that this probability increases with $\alpha$,  but a large $\alpha$ also reduces the stepsize which potentially slows down the convergence. Therefore, the choice of $\alpha$ is very important for the algorithm implementation.
\hfill \qedsymbol



\end{document}